\newcommand{\xba}{\alpha}
\newcommand{\xbb}{\beta}

\newcommand{\xbe}{\in}
\newcommand{\xbf}{\phi}
\newcommand{\xbg}{\gamma}

\newcommand{\xbm}{\mu}

\newcommand{\xbo}{\omega}

\newcommand{\xCN}{\neg}

\newcommand{\xCQ}{\emptyset}

\newcommand{\xCf}{\hspace{0.1em}}

\newcommand{\xcA}{\forall}

\newcommand{\xcN}{\hspace{0.2em}\not\sim\hspace{-0.9em}\mid\hspace{0.8em}}

\newcommand{\xcS}{\bigcap}
\newcommand{\xcT}{\bot}

\newcommand{\xcV}{\bigcup}

\newcommand{\xcb}{\subset}
\newcommand{\xcc}{\subseteq}

\newcommand{\xce}{\not\in}

\newcommand{\xcl}{\vdash}

\newcommand{\xcn}{\hspace{0.2em}\sim\hspace{-0.9em}\mid\hspace{0.58em}}

\newcommand{\xcp}{\rightarrow}

\newcommand{\xcs}{\cap}
\newcommand{\xcu}{\wedge}
\newcommand{\xcv}{\cup}

\newcommand{\xcx}{\Diamond}

\newcommand{\xDH}{\item }

\newcommand{\xDf}{\hspace{0.5em}}

\newcommand{\xda}{{\cal A}}

\newcommand{\xdr}{{\cal R}}

\newcommand{\xEI}{\begin{itemize}}
\newcommand{\xEJ}{\end{itemize}}

\newcommand{\xEP}{ \\ }

\newcommand{\xEd}{\neq}

\newcommand{\xEh}{\begin{enumerate}}

\newcommand{\xEj}{\end{enumerate}}

\newcommand{\xeb}{\prec}

\newcommand{\xej}{\lhd}

\newcommand{\xer}{\sqsubset}

\newcommand{\xfI}{\mbox{I}}

\newcommand{\Xl}{\ldots}

\newcommand{\bl}{\begin{lemma} \rm}
\newcommand{\el}{\end{lemma}}
\newcommand{\br}{\begin{remark} \rm}
\newcommand{\er}{\end{remark}}
\newcommand{\be}{\begin{example} \rm}
\newcommand{\ee}{\end{example}}
\newcommand{\bco}{\begin{corollary} \rm}
\newcommand{\eco}{\end{corollary}}
\newcommand{\bc}{\begin{claim} \rm}
\newcommand{\ec}{\end{claim}}
\newcommand{\bfa}{\begin{fact} \rm}
\newcommand{\efa}{\end{fact}}
\newcommand{\bp}{\begin{proposition} \rm}
\newcommand{\ep}{\end{proposition}}
\newcommand{\bd}{\begin{definition} \rm}
\newcommand{\ed}{\end{definition}}
\newcommand{\bcs}{\begin{construction} \rm}
\newcommand{\ecs}{\end{construction}}
\newcommand{\bcd}{\begin{condition} \rm}
\newcommand{\ecd}{\end{condition}}
\newcommand{\bt}{\begin{theorem} \rm}
\newcommand{\et}{\end{theorem}}
\newcommand{\bn}{\begin{notation} \rm}
\newcommand{\en}{\end{notation}}
\newcommand{\bfi}{\begin{bild} \rm}
\newcommand{\efi}{\end{bild}}
\newcommand{\bsta}{\begin{statement} \rm}
\newcommand{\esta}{\end{statement}}
\newcommand{\bcom}{\begin{comment} \rm}
\newcommand{\ecom}{\end{comment}}
\newcommand{\bdia}{\begin{diagram} \rm}
\newcommand{\edia}{\end{diagram}}

\newcommand{\bfc}{\begin{figure}[htb] \begin{center}}
\newcommand{\efc}{\end{center} \end{figure}}

\documentclass{article}

\usepackage{amssymb,latexsym,epic,eepic,rotating}

\title{A Comment on Argumentation
}

\author{Karl Schlechta
\thanks{
schcsg@gmail.com - https://sites.google.com/site/schlechtakarl/ -
Koppeweg 24, D-97833 Frammersbach, Germany}
\thanks{
Retired, formerly: Aix-Marseille Universit\'{e}, CNRS, LIF UMR 7279, F-13000
Marseille, France
}
}

\begin{document}

\newtheorem{lemma}{Lemma}[section]
\newtheorem{theorem}[lemma]{Theorem}
\newtheorem{proposition}[lemma]{Proposition}
\newtheorem{corollary}[lemma]{Corollary}
\newtheorem{claim}[lemma]{Claim}
\newtheorem{fact}[lemma]{Fact}
\newtheorem{remark}[lemma]{Remark}
\newtheorem{definition}{Definition}[section]
\newtheorem{construction}{Construction}[section]
\newtheorem{condition}{Condition}[section]
\newtheorem{example}{Example}[section]
\newtheorem{notation}{Notation}[section]
\newtheorem{bild}{Figure}[section]
\newtheorem{comment}{Comment}[section]
\newtheorem{statement}{Statement}[section]
\newtheorem{diagram}{Diagram}[section]

\renewcommand{\labelenumi}
  {(\arabic{enumi})}
\renewcommand{\labelenumii}
  {(\arabic{enumi}.\arabic{enumii})}
\renewcommand{\labelenumiii}
  {(\arabic{enumi}.\arabic{enumii}.\arabic{enumiii})}
\renewcommand{\labelenumiv}
  {(\arabic{enumi}.\arabic{enumii}.\arabic{enumiii}.\arabic{enumiv})}

\maketitle

\setcounter{secnumdepth}{3}
\setcounter{tocdepth}{3}

\begin{abstract}

We use the theory of defaults and their meaning of
\cite{GS16} to develop (the outline of a) new theory of argumentation.

\end{abstract}

\tableofcontents

%
%
%
\section{
Introduction
}
\subsection{
Abstract description
}

Argumentation is about putting certain objects together.

There are three things to consider:
 \xEh
 \xDH the objects themselves, and their inner structure (if they have any)
-
this inner structure may be revealed successively, or be immediately
present,
 \xDH rules about how to put them together,
 \xDH avoid certain results (contradictions) in the resulting pattern.
 \xEj

To help intuition, we picture as result of an argumentation, an
inheritance
network the agents can agree on.

This network may consist of strict
and defeasible rules only, with no elements or sets it is applied to.
Think of the argumentation going on when writing a book about
medical diagnosis. This will not be about particular cases, but about
strict and default rules. ``Sympton $x$ is usually a sign of illness $y,$
but there are the following exceptions:  \Xl''
In addition, the network might contain cycles. There is nothing wrong with
cycles.
Mathematics is full of cycles, equivalences and their proofs. But consider
also the following: We work in the set of adult land mammals. ``Most
elefants
weigh more than 1 ton.'' ``Most elements (i.e. adult land mammals) which
weigh more than 1 ton are elefants.'' There is nothing in principle wrong
with
this either - except, in reality, we forgot perhaps about hippopotamus
etc.

Arguments need not be contradictions to what exists already. They can be
confirmations, elaborations, etc. For instance, we might have the default
rule that birds fly, and clarify that penguins don't fly. This is not a
contradiction, but an elaboration.
\subsection{
The structure of the objects
}

Facts are either so simple that a dispute seems unreasonable. Or, they are
a combination of basic facts and (default) rules, like, what $ \xfI $ see
through
my microscope is really there, and not an artifact of some speck of
dust on the lenses. For simplicity, facts will be basic, undisputable
facts.

Expert opinion may be considered a default rule, where details stay
unexplained, perhaps even unexplainable by the expert himself.

Rules (classical or defaults) have three aspects:
 \xEh
 \xDH the rule itself,
 \xDH the application of the rule,
 \xDH the result of the application of the rule.
 \xEj

Classical rules cannot be contested. We can contest their application,
i.e. one of their prerequisites, or their result, and, consequently,
their application. We can confirm their result by different means,
likewise, their application.

Default rules are much more complicated, but not fundamentally different.
Again, we can attack their
application, by showing that one of the prerequisites does not hold, or,
that we are in an (known) exceptional case. We can attack the conclusion,
and, consequently, the rule, or its application. In particular, we may
attack the conclusion, without attacking the application or the rule
itself,
by arguing that we are in a surprising exceptional case - and perhaps try
to find a new set of exceptions. We can attack the default rule itself,
as in the case of ``normally, tigers are vegans''. We can confirm a rule
by confirming its conclusion, or adding a new rule, which gives the same
result. We can elaborate a default rule, by adding an exception set,
stating
that all exceptions are known, and give the list of exceptions, etc.
We can stop homogenousness (downward inheritance) e.g. for Quakers which
are Republicans, we stop inheriting pacifism (or its opposite).
This is not a contradiction to the default itself, but to the
downward inheritance of the default (or to homogenousness) by
meta-default,
to be precise.

Obviously, the more we add (possible) properties to the objects (here
default
rules), the more we
can attack, elaborate, confirm.

In the following section, we describe our general picture:
 \xEh
 \xDH there is an
arbiter which checks for consistency, and directs the discussion,
 \xDH how to
handle classical arguments and resulting contradictions,
 \xDH how to handle default arguments.
 \xEj
\section{
The classical part
}

We suppose there is an arbiter, whose role is to check consistency, and to
authorise participants to speak.

If the arbiter detects an inconsistency, then he points out the
``culprits'',
i.e. minimal inconsistent sets. As he detects inconsistencies immediately,
the last argument will be in all those sets. The last argument need not
be the problem, it might be one of the earlier arguments.

He asks all participants if they wish to retract one of the arguments
involved in at least one minimal inconsistent set.
(They have to agree unanimously on such retraction.)
If there is no minimal
inconsistent set left, the argumentation proceeds with the ``cleaned''
set of arguments, as if the inconsistency did not arise.
Of course, arguments which were based on some of the retracted arguments
are now left ``hanging in the air'', and may be open to new attacks.

If not, i.e. at least one minimal inconsistent set is left, the
participants
can defend (and attack) the arguments involved in those sets.
The arbiter will chose the argument to be attacked/defended.
See Example 
\ref{Example Symmetrical} (page 
\pageref{Example Symmetrical})  below.
Suppose $ \xba $ is the argument chosen, then a defense will try to prove
or
argue for $ \xba,$ an attack will try to prove or argue for $ \xcx \xCN
\xba,$ i.e.
it is possible or consistent that $ \xCN \xba.$
In particular, an attacker might try to prove $ \xcT,$ or some other
unlikely
consequence of $ \xba $ (and some incontested $ \xbb $'s), and he need
not begin
with some $ \xba \xcp \xbg,$ it might be a more roundabout attack.

If at least one minimally inconsistent set is left with all elements
defended, then there is a deadlock, and the arbiter declares failure.

Consider

\be

$\hspace{0.01em}$


\label{Example Symmetrical}

We argue semantically. Let $A:=\{x,a\},$ $B:=\{x,b\},$ $C:=\{x,c\},$
$Y:=\{a,b,c\}.$
Let $Y$ be the last set added. For $A,B,C,$ the situation is symmetrical.
Let $Z \xEd Z' $ be $A,B,$ or $C,$ then $Y \xcs Z \xEd \xCQ,$ but $Y \xcs
Z \xcs Z' = \xCQ,$ $Z \xcs Z' =\{x\},$ etc.
Moreover, $A \xcs B \xcc C,$ etc. Thus, $ \xCf A$ and $B$ together are an
argument for $C,$ etc.,
so they argue for each other, and there is no natural way to chose any
of $A,B,C$ to be attacked. Thus, it is at the discretion of the parties
involved (or the arbiter) to chose the aim of any attack - apart from $Y,$
which is not supported by any of $A,B,C.$ Still, $Y$ might in the end be
the strongest argument.

We may add $D,$ $E,$ with $D:=\{x,d\},$ $Y:=\{a,b,c,d\}$ etc., the example
may be extended
to arbitrarily many sets.

\ee

At any moment, any argument can be attacked, not only if an inconsistency
arises. We may continue an argumentation, even if not all minimally
inconsistent subsets are treated as yet, but the arbiter has to keep track
of
them, and of the use of their elements. They and their consequences may
still be questioned.
\section{
Defaults
}
\subsection{
The classical part of defaults
}

We see defaults primarily not as rules, but as relatively complicated
classical constructions, which we may see as objects for the moment.
The default character is in applying those objects, not in the
objects themselves.

We follow here the theory described in Chapter 11 of
 \cite{GS16}, and
for convenience of the reader, we repeat in
Section \ref{Section Appendix} (page \pageref{Section Appendix})  (essentially)
Section 11.4 of  \cite{GS16}.

In our view, a (semantical) default $(\xCf X:Y)$ says:

 \xEh
 \xDH ``most'' elements of $X$ are in $Y,$
 \xDH there may be exception sets $X_{1},$ $X_{2},$ etc. of $X,$ where the
elements
are ``mostly'' not in $Y$ (but $X_{1} \xcv X_{2} \xcv  \Xl.$ has to be a
``small'' subset of $X),$
 \xDH in addition, there may be a ``very small'' subset $X' \xcc X,$ which
contains
``surprise elements'' (i.e. not previously known exceptions), which are not
in $Y,$
 \xDH in addition, we may require that subsets of $X$ ``normally''
behave in a homogenous way.
 \xEj

The notions of ``most'', ``small'' etc. are left open, a numerical
interpretation
suffices for the intuition. These notions are discussed in depth e.g. in
 \cite{GS08f} and  \cite{GS10}.

Introducing a default has to result in a (classically) consistent
theory. E.g., it must not be the case that $ \xcA x \xbe X.x \xce Y,$ this
contradicts
the first requirement about defaults (and any reasonable interpretation
of ``most'').
\subsection{
The default part of defaults
}

This leads to a hierarchy as defined in
Section 
\ref{Section Form-Sol-V2} (page 
\pageref{Section Form-Sol-V2})  below. We use the hierarchy to define
the $ \xCf use$
of the defaults.

To use the standard example with birds, penguins, fly, we proceed as
follows. Suppose we introduce a bird $x$ into the discussion.
We try to put $x$ as low as possible in the hierarchy, i.e. into the
set of birds, but not into any known exception set, and much less into
any ``surprise'' set. Only (classical) inconsistency, as checked by the
arbiter, may force us to climb higher.
Thus, unless there is a contradiction, we let $x$ fly.
\subsection{
Attacks against defaults and their conclusions
}

Classical rules are supposed to be always true. Thus, classical
rules themselves cannot be attacked, and an attack against a
classical conclusion has to be an attack against one of its
prerequisites.

Attacks against defaults can be attacks against
 \xEh
 \xDH the rule itself,
 \xDH one of the prerequisites,
 \xDH membership in or not in one of the exception sets,
 \xDH membership in or not in the surprise set,
 \xDH perhaps even the notions of size involved,
 \xDH etc.
 \xEj

Each component of a default rule may be attacked.
\section{
Comments
}

We assume that there
is no fundamental difference between facts and conclusions:
Usually, we were told facts, remember facts, have read facts, observed
facts (perhaps with the help of a telescope etc.). These things
can go wrong. Situations where things are obvious, and no error
seems humanly possible, will not be contradicted.
\subsection{
Auxiliary elements
}

We now introduce some auxiliary elements which may help in the
argumentation.

 \xEh
 \xDH ``$ \xfI $ agree.''

This makes an error in this aspect less likely, as both parties agree -
but
still possible!

 \xDH ``$ \xfI $ confirm.''

$ \xfI $ am very certain about this aspect.

 \xDH
Expert knowledge:

Expert knowledge and its conclusions act as ``black box defaults'', which
the expert himself may be unable to analyse. Other experts (in the same
field)
will share the conclusion. (This is simplified, of course.)

(One way to contest an expert's
conclusion is to point out that he neglected an aspect of the situation,
which is outside his expertise. His ``language of reasoning'' is too poor
for the situation.)

 \xDH
The arbiter may ask questions.

 \xEj
\subsection{
Examples of attacks
}

 \xEh
 \xDH
Defaults:

Normally, there is a bus line number 1 running every 10 minutes between 10
and
11 in the morning.

Attack: No, the conclusion is wrong.

Question: Why?

Elaboration:

 \xEh
 \xDH
No, the default is wrong (e.g.: it is line number 2 running every 10
minutes).
 \xDH
Yes, but this is not homogenous, i.e. does not break down to subsets,
and we know more. (For instance, we know that today is Tuesday or
Wednesday,
and it runs that often only Monday, Thursday, Friday, Saturday, Sunday -
but
we do not know this, only that is does not apply to all days of the week.)

 \xDH
Yes, but today is an exception, and we know this. (e.g., we know that
today is Tuesday, and we know that Tuesday is an exception.)
(In addition, there might be exceptional Tuesdays, Christmas market day,
etc.  \Xl)

 \xDH
Yes, but $ \xfI $ do not know why this is an exception. (This is a
surprise case,
$ \xfI $ know about different days, but today should not be an exception,
still
$ \xfI $ was just informed that it does not hold today.) We do not attack
the
default, nor the applicability - but agree that it fails here.

 \xEj

 \xDH
Classical conclusions:

From $ \xCf A$ and $B,$ $C$ follows classically.

Attack: $C$ does not hold.

Question: Why?

Elaboration:

 \xEh
 \xDH
$ \xCf A$ does not hold or $B$ does not hold, but $ \xfI $ do not know
which.

 \xDH
$ \xCf A$ does not hold.

 \xDH
$B$ does not hold.

 \xDH
$ \xCf A$ does not hold, and $B$ does not hold.

 \xEj

 \xDH
Fact: $ \xCf A$ holds.

Attacks: No, $ \xCf A$ does not hold.

Question: Why?

Elaboration:

 \xEh
 \xDH
You remember incorrectly.

 \xDH
You did not observe well.

 \xDH
Your observation tools do not work.
 \xDH
You were told something wrong.
 \xDH
etc.
 \xEj

 \xDH
Expert knowledge, expert concludes that $ \xCf A.$

Attack: $ \xCf A$ does not hold.

Question: Why?

Elaboration:

The situation involves aspects where you are not an expert. It is beyond
your language. (Of course, the expert can ask for elaboration  \Xl.)

 \xEj
\clearpage
\section{
Appendix - Section 11.4 in \cite{GS16}
}

\label{Section Appendix}

\bd

$\hspace{0.01em}$


\label{Definition Attached}

If $ \xba \xcn \xbf $ or $ \xba \xcN \xbf,$ we say that the default $
\xbf $ is attached \index{attached}  to $ \xba.$

Given any fixed default theory, let $ \xda $ be the set of $ \xba $ (or
$M(\xba)),$
to which some default is attached

\ed

We work here in propositional logic, and on the semantic level.

We assume a classical background theory $B,$ and a set of
classical formulas $ \xba_{1}, \Xl, \xba_{n}$ to which defaults are
attached,
$ \xba_{i} \xcn \xbf_{i,1}, \Xl, \xba_{i} \xcn \xbf_{i,j_{i}},$ where
some or all of the $ \xcn $ may also be $ \xcN $
(without being necessarily $ \xcn \xCN).$

\bcd

$\hspace{0.01em}$


\label{Condition Consistency}

We assume the following consistency conditions: \index{consistency conditions}
 \xEh
 \xDH $B$ is classically consistent
 \xDH For each $ \xba_{i}$ the defaults attached to $ \xba_{i}$ together
with $B$ are
jointly consistent.

In particular, the theory of defaults attached to one of the $ \xba_{i}$
must be
consistent.
For instance, $ \xba_{i} \xcl \xba_{j}$ and $ \xba_{i} \xcn \xCN \xba_{j}$
together are inconsistent.
We thus rule out default theories like $\{: \xbf / \xbf,$ $: \xCN \xbf /
\xCN \xbf \}.$

A negated default like $ \xba \xcN \xbf $ needs a model of $ \xba \xcu
\xCN \xbf $ to be consistent,
so we replace $ \xba \xcN \xbf $ by $ \xba \xcu \xCN \xbf $ for the
consistency check.
 \xEj

\ecd

The $ \xcl $-relation between the $ \xba_{i}$'s gives a specificity relation
\index{specificity relation}  by strict
inclusion. We use this for the inheritance relation, and to solve (some)
conflicts.

\bd

$\hspace{0.01em}$


\label{Definition Inheritance}

This definition describes how to obtain a consistent default theory
at every point in the universe, using a theory revision approach, with
specificity solving some conflicts.
Of course, modifications are possible, but the general idea seems
sound.

We define the set of valid defaults at some point.
This influences the relation $ \xer $ as defined in
Definition 
\ref{Definition Preference} (page 
\pageref{Definition Preference}), but not the relation $ \xej $ as
defined in
Definition 
\ref{Definition RelevantOrder-V2} (page 
\pageref{Definition RelevantOrder-V2}),
as the latter depends only on
the sets to which defaults are attached, and not which defaults
are attached.

We consider now some classical formula $ \xbb $ - it may be one of the $
\xba_{i}$'s
to which defaults are attached, or not.
 \xEh
 \xDH Visible defaults at $ \xbb $
 \index{visible defaults}

The defaults visible at $ \xbb:$ All defaults attached to some $
\xba_{i}$ such that
$ \xbb \xcl \xba_{i}$ are considered visible at $ \xbb.$

 \xDH Valid defaults at $ \xbb $
 \index{valid defaults}
 \xEh
 \xDH The visible defaults are ordered by the $ \xcl $ relation between
the $ \xba $'s to
which they are attached.
The more specific
ones are considered stronger defaults (for this $ \xbb).$

Of course, we can plug in here any partial relation, if it seems more
suitable, e.g. $ \xcn $ itself, as is often done in defeasible inheritance
networks.
 \xDH Consider now the set of visible defaults, together with the
classical
information available at $ \xbb.$ (If $ \xbb' \xcl \xbb,$ but not
conversely, we need not
consider $ \xbb',$ etc.) If this set is inconsistent:
 \xEh
 \xDH Consider first in parallel all minimal inconsistent subsets
involving
classical information.
They must contain at least one default, as the classical information
was supposed to be consistent.
Eliminate simultaneously from each such set the weakest (by the $ \xcl $
relation)
defaults. (If
there are several weakest ones, eliminate all the weakest ones.)

As contradictions involving classical, i.e., strongest information, seem
to
be more serious, we do these sets first.

 \xDH Consider now all remaining minimal inconsistent subsets of default
information only. (Note that some might already have been eliminated in
the previous step, but for other reasons.)

Proceed as in the previous step, i.e. eliminate the weakest information.

 \xDH We call the remaining defaults, visible at $ \xbb,$ the defaults
valid at $ \xbb.$

We will now work with the valid defaults only.
 \xEj
 \xEj
 \xEj

\ed

We now define ``relevant'' sets, sets where some default may change.
We work in some fixed default theory, recall from
Definition 
\ref{Definition Attached} (page 
\pageref{Definition Attached})  that $ \xda $ is the set of $ \xba $
or $M(\xba)$ to which
at least one default is attached.

\bd

$\hspace{0.01em}$


\label{Definition Relevant}

Consider the model variant of $ \xda.$

Let $S(\xda)$ \index{$S(\xda)$}  $:=$ $\{ \xcS \xda':$ $ \xda' \xcc
\xda \}$
and $U(\xda)$ \index{$U(\xda)$}  $:=$ $\{ \xcV \xda':$ $ \xda' \xcc
\xda \},$

Let $ \xdr:=\{X- \xCf Y:$ $X-Y \xEd \xCQ,$ $X \xbe S(\xda),$ $Y \xbe
U(\xda)\},$

where $X=U$ (the universe), and $Y= \xCQ $ are possible.

$R \xbe \xdr $ is called a relevant set (or formula).

Let $ \xdr_{f}$ \index{$ \xdr_{f}$}  be the $ \xcc $-minimal elements of $
\xdr $ - ``f'' for finest.

\ed

Obviously, the elements of $ \xdr_{f}$ are pairwise disjoint.
The motivation for $ \xdr_{f}$ is the following. Let $X-Y \xbe \xdr_{f},$
then
all elements of $X-Y$ should satisfy all defaults valid at $X,$
but no other, as they are in no other sets to which defaults are
attached.

\be

$\hspace{0.01em}$


\label{Example Relation-V1}

Let $A,A' \xcc U,$ $A'' \xcc A',$ let defaults be attached to $A,A',A''
,$
let all resulting intersections and set-differences
be non-empty, if possible, e.g. $A'' -A \xEd \xCQ.$

See Diagram 
\ref{Diagram Form-Sets} (page 
\pageref{Diagram Form-Sets})  for illustration.

So $ \xda =\{A,A',A'' \},$

$ \xdr =\{U,$ $U- \xCf A,$ $U-A',$ $U-A'',$ $U-A-A',$

$ \xCf A,$ $A-A',$ $A-A'',$

$A',$ $A' - \xCf A,$ $A' -A'',$ $A' -A-A'',$

$A'',$ $A'' - \xCf A,$

$A \xcs A',$ $A \xcs A' -A'',$

$A \xcs A'' \}$

and

$ \xdr_{f}=\{U-A-A',$

$A-A',$

$A' -A-A'',$

$A'' - \xCf A,$

$A \xcs A' -A'',$

$A \xcs A'' \}$

It is useful to code the sets $ \xCf A,$ $A',$ $A'' $ by 3 bits, e.g. the
left bit codes $ \xCf A,$ the middle one $A',$ the right one $A''.$ We
then have

$ \xdr_{f}=\{U-A-A' =000,$ $A-A' =100,$ $A \xcs A' -A'' =110,$ $A \xcs A''
=111,$ $A'' -A=011,$ $A' -A-A'' =010\},$

where the codes 001, 101 do not exist, as $A'' \xcc A'.$

\ee

We will order the models by groups, and within groups by their ``quality''
satisfying defaults. We first define the latter relation, to be
denoted $ \xer.$ It will be used to order the models within the sets $
\xbo (A),$
see Construction 
\ref{Construction Preference-V2} (page 
\pageref{Construction Preference-V2}).

\bd

$\hspace{0.01em}$


\label{Definition Preference}

Suppose we are at some relevant $ \xbb,$ with the valid defaults $
\xbf_{1}, \Xl, \xbf_{n}.$
We order the models of $ \xbb $ according to satisfaction of the $
\xbf_{i}.$
There are different possibilities to define $ \xer $: \index{$ \xer $}

 \xEh
 \xDH By subsets: if the set of $ \xbf_{i}$ satisfied by $m$ is a subset
of those satisfied
by $m',$ then $m' $ is better than $m,$ $m' \xer m.$

 \xDH By cardinality: if the set of $ \xbf_{i}$ satisfied by $m$ is a
smaller
than the set of those satisfied
by $m',$ then $m' $ is better than $m,$ $m' \xer m.$

 \xDH Some more complicated order, preferring certain defaults over
others.
(This might be interesting for contrary-to-duty obligations. Note that
we solve here the case of additional information in real-life
situations, too.)
 \xDH In particular, we may order the valid defaults by specificity
of the $ \xba $ to which they are attached,
and satisfy the most specific ones first, then the next, etc.,
resulting in a lexicographic order.
 \xEj

\ed

Note that the construction is more adapted to sets with
rare exceptions, than
to a classification, like vertebrates into fish, amphibians, reptiles,
birds,
and mammals, where all are exceptions - and those outside these subsets
are
the real exceptions.



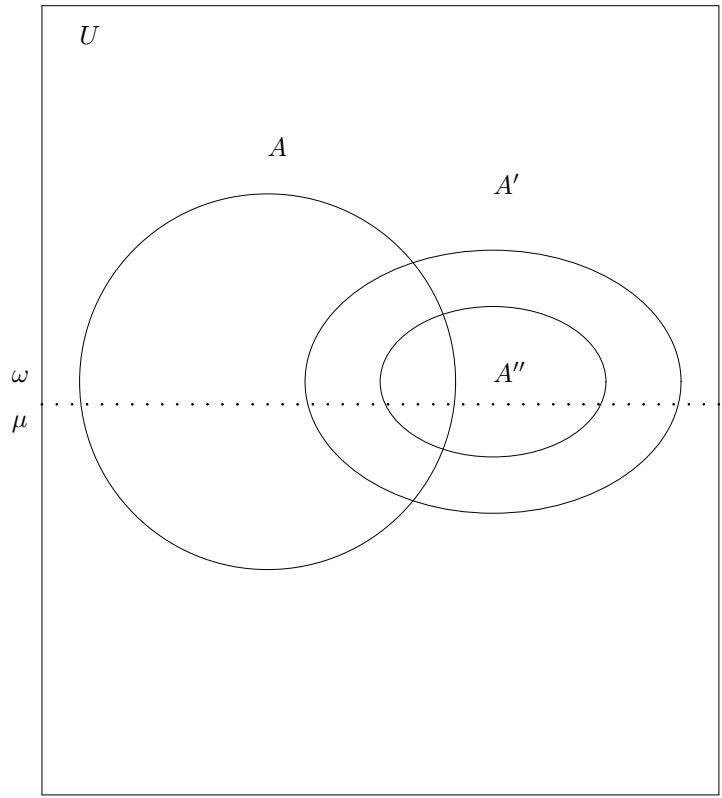
\vspace{30mm}

\begin{diagram}

\label{Diagram Form-Sets}
\index{Diagram Example Form Sets}

\centering
\setlength{\unitlength}{1mm}
{\renewcommand{\dashlinestretch}{30}
\begin{picture}(110,150)(0,0)

\path(10,25)(10,130)
\path(100,25)(100,130)
\path(10,25)(100,25)
\path(10,130)(100,130)

\multiput(10,77)(2,0){46}{\circle*{0.2}}

\put(40,80){\circle{50}}

\put(70,80){\ellipse{50}{35}}
\put(70,80){\ellipse{30}{20}}

\put(15,125){$U$}
\put(40,110){$A$}
\put(70,105){$A'$}
\put(70,80){$A''$}

\put(6,74){$\xbm$}
\put(6,80){$\xbo$}

\put(10,15){The dotted line separates $\xbo (X)$ from $\xbm (X)$}

\put(40,8){{\rm\bf Example, the sets}}

\end{picture}
}

\end{diagram}

\vspace{4mm}

\vspace{3mm}


\vspace{3mm}

\vspace{3mm}


\vspace{3mm}


\vspace{3mm}


\vspace{3mm}


\vspace{3mm}


\vspace{3mm}



\vspace{3mm}


\vspace{3mm}

\subsection{
The construction
}

\label{Section Form-Sol-V2}

Note that, again, the following definition is independent of the actual
defaults,
and uses only the fact that defaults are attached to the $A \xbe \xda.$
The order to be constructed, $ \xej,$ will be used as scaffolding in the
construction of $ \xeb $ in Construction 
\ref{Construction Preference-V2} (page 
\pageref{Construction Preference-V2}),
which orders
packets of models.

(See Example \ref{Example Relation-V2} (page \pageref{Example Relation-V2})  and
Diagram \ref{Diagram Form-Sets} (page \pageref{Diagram Form-Sets})  and
Diagram \ref{Diagram Set-Rel-V2} (page \pageref{Diagram Set-Rel-V2}).)

Recall Definition 
\ref{Definition Relevant} (page 
\pageref{Definition Relevant}).

\bd

$\hspace{0.01em}$


\label{Definition RelevantOrder-V2}

We order $ \xdr_{f}$ \index{$ \xdr_{f}$}  as follows:

$X \xej Y$ \index{$X \xej Y$}  iff $\{A_{i}:X \xcc A_{i}\} \xcb \{A_{i}:Y
\xcc A_{i}\}$ for $X,Y \xbe \xdr_{f},$

i.e. by the subset relation, of the $A_{i}$ they are in.

\ed

This order expresses exceptionality of the $X \xbe \xdr_{f}.$ If $X \xej
Y,$ then
$Y$ is a subset of more (by the subset relation) $A_{i}$ than $X$ is,
defaults were attached to the $A_{i},$ so the $A_{i}$ are sets of
exceptions.
For instance, $A_{0}$ might be the set of birds, $A_{1}$ the set of
penguins,
being a penguin (an element of $A_{0}$ and $A_{1})$
is more exceptional than being a bird (an element of $A_{0}$ only).

Recall the coding of the
$X \xbe \xdr_{f}$ in Example 
\ref{Example Relation-V1} (page 
\pageref{Example Relation-V1}).
This example is continued
in Example \ref{Example Relation-V2} (page \pageref{Example Relation-V2}), and
$ \xej $ is then (the transitive closure of)
$U-A-A' -A'' $ (000) $ \xej $ $A-A' $ (100) $ \xej $ $A \xcs A' -A'' $
(110) $ \xej $ $A \xcs A'' $ (111),
$U-A-A' -A'' $ (000) $ \xej $ $A' -A-A'' $ (010) $ \xej $ $A'' -A$ (011) $
\xej $ $A \xcs A'' $ (111),
$A' -A-A'' $ (010) $ \xej $ $A \xcs A' -A'' $ (110).

See Diagram \ref{Diagram Set-Rel-V2} (page \pageref{Diagram Set-Rel-V2}).

We are now ready to construct the preferential relation $ \xeb $ between
model sets, it is implicitely extended to their elements.

Take e.g. in Example 
\ref{Example Relation-V2} (page 
\pageref{Example Relation-V2})  the set $A \xcs A' -A''
.$ It will ``see'' the
defaults valid for $A \xcs A',$ but not necessarily those for $A''.$
Every element of $A \xcs A' -A'' $ can satisfy all those defaults, or
only
a part of them (or none). The best elements are those which satisfy all,
the worst those which satisfy none. The precise relation is described by
$ \xer,$ see Definition 
\ref{Definition Preference} (page 
\pageref{Definition Preference}). Let us call $ \xbm
(X)$ the set of those
elements which satisfy all defaults, the set of the others $ \xbo (X),$
for given $X.$
Elements which do not satisfy all defaults are doubly exceptional, as
they are not in any subclass where this failure is ``explained''.
They are the ``unexcused'' exceptions,
they have no reason to be exceptional, they are surprises,
and doubly exceptional.
So they should sit higher up in the hierarchy.

We do not think that there is a unique reasonable solution.
Two ideas come to mind:
 \xEh
 \xDH
Put all sets $ \xbo (X)$ above all other
elements, in their own hierarchy, defined by $ \xej,$ and ordered
internally
by $ \xer.$ This is the radical approach.
 \xDH
A less radical idea is to put them above the immediate $ \xej $-successors
of $X.$
E.g. $ \xbo (A' -A-A'')$ will sit above $A \xcs A' -A'' $ and $A' \xcs
A'' -A$ - or, rather,
above $ \xbm (A \xcs A' -A'')$ and $ \xbm (A' \xcs A'' - \xCf A).$
This is what we will do here.
 \xEj

\bcs

$\hspace{0.01em}$


\label{Construction Preference-V2}

(See Example \ref{Example Relation-V2} (page \pageref{Example Relation-V2})  and
Diagram \ref{Diagram Mod-Rel-V2} (page \pageref{Diagram Mod-Rel-V2}).)

We define a ``packetwise'' order $ \xeb $ \index{$ \xeb $}  between model
sets, it will then be
extended to the elements.

The relation $ \xeb $ has two parts:
 \xEh
 \xDH
$ \xbm (X) \xeb \xbm (Y)$ iff $X \xej Y$

for $X,Y \xbe \xdr_{f}$

 \xDH
$ \xbm (X) \xeb \xbo (Y)$ iff $X \xej Y$ or $X=Y$ or $X$ is a direct $
\xej $-successor of $Y$

for $X,Y \xbe \xdr_{f}.$

 \xEj

Note:
 \xEh
 \xDH
We do not continue the order above the $m \xbe \xbo (A),$ they are maximal
elements -
except for the interior order $ \xer $ among themselves.

 \xDH
The $ \xbm (X)$ are ``flat'', all are best possible,
There is no interior order in this set $ \xbm (X).$

 \xDH
The inside structure
of the $ \xbo (X)$ is given by $ \xer $
according to Definition 
\ref{Definition Preference} (page 
\pageref{Definition Preference})  for this $X.$

 \xDH
If any of the sets is empty, we just
omit it and close under transitivity.
 \xEj

\ecs

More formally, for $X,X' \xbe \xdr_{f},$ and element wise:
\begin{flushleft}
\[ m \xeb m':= \left\{ \begin{array}{ll}
\mbox{there is} \xDf X \xDf \mbox{and} \xDf m,m' \xbe \xbo (X)
\xDf \mbox{and} \xDf m \xer m' \xEP
or \xEP
\mbox{there are} \xDf X,X' \xDf \mbox{and} \xDf m \xbe \xbm (X)
\xDf \mbox{and} \xDf m' \xbe \xbm (X') \xDf
\mbox{and} \xDf X \xej X' \xEP
or \xEP
\mbox{there are} \xDf X,X' \xDf \mbox{and} \xDf m \xbe
\xbm (X) \xDf \mbox{and}
\xDf m' \xbe \xbo (X') \xDf \mbox{and} \xEP
(\xDf X \xej X' \xDf \mbox{or} \xDf X=X' \xDf \mbox{or} \xDf
X \xDf \mbox{is a direct} \xej- \mbox{successor}
\xDf \mbox{of} \xDf X') \xEP
\end{array}
\right.
\]
\end{flushleft}

 \index{$ \xeb $}

The final order between models is thus basically lexicographic with two
parts,
first by $ \xej,$ then by $ \xer.$

\be

$\hspace{0.01em}$


\label{Example Relation-V2}

We continue
Example \ref{Example Relation-V1} (page \pageref{Example Relation-V1}),
see also Diagram \ref{Diagram Form-Sets} (page \pageref{Diagram Form-Sets}).

Recall that $ \xda =\{A,A',A'' \},$ and
$ \xdr_{f}=\{U-A-A' -A'' =000,$ $A-A' =100,$ $A \xcs A' -A'' =110,$ $A
\xcs A'' =111,$ $A'' -A=011,$ $A' -A-A'' =010\}.$
 \xEh
 \xDH
We define $ \xej $ by
Definition 
\ref{Definition RelevantOrder-V2} (page 
\pageref{Definition RelevantOrder-V2}):
$U-A-A' -A'' =000$ $ \xej $ $A-A' =100$ $ \xej $ $A \xcs A' -A'' =110$ $
\xej $ $A \xcs A'' =111,$
$U-A-A' -A'' =000$ $ \xej $ $A' -A-A'' =010$ $ \xej $ $A'' -A=011$ $ \xej
$ $A \xcs A'' =111,$
$A' -A-A'' =010$ $ \xej $ $A \xcs A' -A'' =110,$
closed by transitivity.

See Diagram \ref{Diagram Set-Rel-V2} (page \pageref{Diagram Set-Rel-V2}).

 \xDH
We construct the relation between models according to
Construction 
\ref{Construction Preference-V2} (page 
\pageref{Construction Preference-V2}).
as follows (considering only the
model sets, and neglecting the internal structure of the $ \xbo (X)).$
Constructing
the full relation (and closing under transitivity) is then trivial.

$ \xbm (000) \xeb \xbm (100) \xeb \xbm (110) \xeb \xbm (111),$
$ \xbm (000) \xeb \xbm (010) \xeb \xbm (011) \xeb \xbm (111),$
$ \xbm (010) \xeb \xbm (110),$

$ \xbm (100) \xeb \xbo (000),$
$ \xbm (010) \xeb \xbo (000),$
$ \xbm (110) \xeb \xbo (100),$
$ \xbm (110) \xeb \xbo (010),$
$ \xbm (011) \xeb \xbo (010),$
$ \xbm (111) \xeb \xbo (110),$
$ \xbm (111) \xeb \xbo (011).$

Finally, we close under transitivity.

See Diagram \ref{Diagram Mod-Rel-V2} (page \pageref{Diagram Mod-Rel-V2}).

 \xEj


\vspace{30mm}

\begin{diagram}

\label{Diagram Set-Rel-V2}
\index{Diagram Example Set Rel V2}

\centering
\setlength{\unitlength}{1mm}
{\renewcommand{\dashlinestretch}{30}
\begin{picture}(110,115)(0,0)

\put(40,15){$U-A-A'-A''=000$}

\put(10,45){$A-A'-A''=100$}
\put(70,45){$A'-A-A''=010$}

\put(40,75){$A \xcs A'-A''=110$}
\put(100,75){$A' \xcs A''-A=011$}

\put(70,105){$A \xcs A' \xcs A''=111$}

\path(52,20)(25,43)
\path(58,20)(82,43)

\path(25,50)(52,73)

\path(82,50)(58,73)
\path(88,50)(115,73)

\path(58,80)(82,103)

\path(115,80)(88,103)

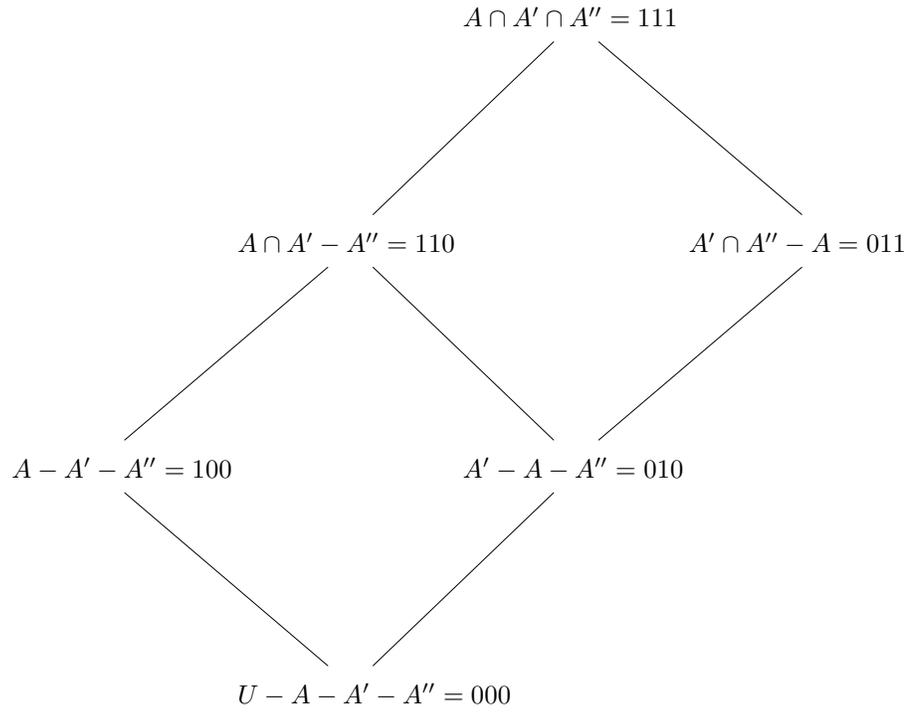
\put(10,5) {{\rm\bf Relation between sets}}


\end{picture}
}

\end{diagram}

\vspace{4mm}

\vspace{3mm}


\vspace{3mm}


\vspace{30mm}

\begin{diagram}

\label{Diagram Mod-Rel-V2}
\index{Diagram Example Mod Rel V2}

\centering
\setlength{\unitlength}{1mm}
{\renewcommand{\dashlinestretch}{30}
\begin{picture}(110,130)(0,0)

\put(50,15){$\xbm (000)$}

\put(20,45){$\xbm (100)$}
\path(23,50)(23,58)
\put(20,60){$\xbo (000)$}

\put(80,45){$\xbm (010)$}
\path(85,50)(85,58)
\put(80,60){$\xbo (000)$}

\put(50,75){$\xbm (110)$}
\path(55,80)(55,88)
\put(50,90){$\xbo (010)$}
\path(52,80)(43,88)
\put(35,90){$\xbo (100)$}

\put(110,75){$\xbm (011)$}
\path(117,80)(117,88)
\put(110,90){$\xbo (010)$}

\put(80,105){$\xbm (111)$}
\path(82,110)(65,118)
\put(60,120){$\xbo (110)$}
\path(88,110)(102,118)
\put(100,120){$\xbo (011)$}
\path(85,110)(85,118)
\put(80,120){$\xbo (111)$}

\path(52,20)(25,43)
\path(58,20)(82,43)

\path(25,50)(52,73)

\path(82,50)(58,73)
\path(88,50)(115,73)

\path(58,80)(82,103)

\path(115,80)(88,103)

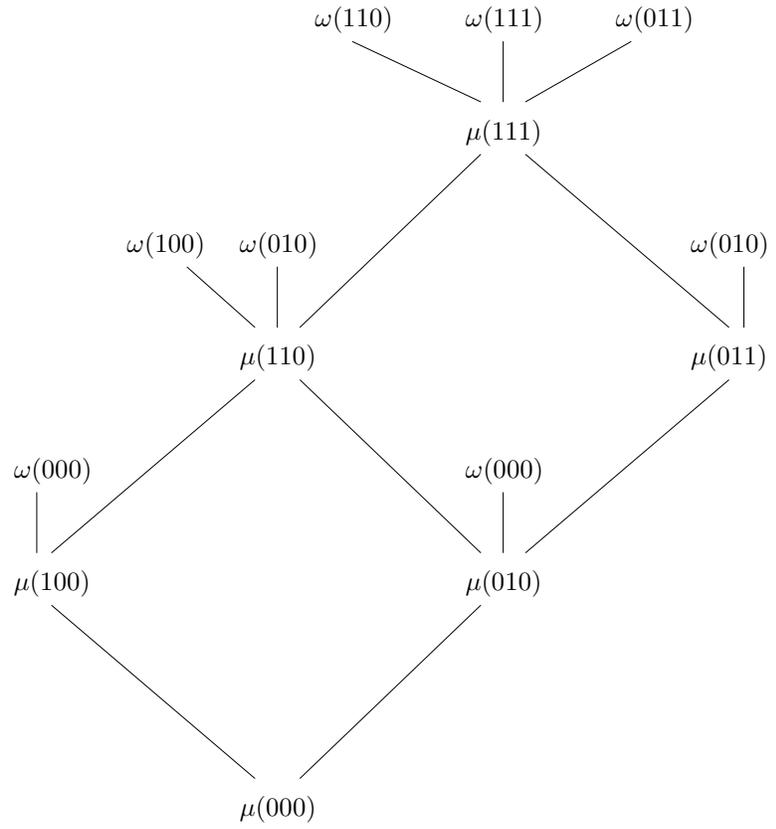
\put(10,5) {{\rm\bf Relation between model sets}}


\end{picture}
}

\end{diagram}

\vspace{4mm}

\vspace{3mm}


\vspace{3mm}

\clearpage

\ee

\end{document}